***Train Timetabling on double track and multiple station capacity railway with useful upper and lower bounds***

Afshin Oroojlooy Jadid[1], Kourosh Eshghi[2]

[1] Department of Industrial and Systems Engineering, Lehigh University

[2]Department of Industrial Engineering, Sharif University of Technology, Tehran, Iran

[1] Corresponding author
Afshin Oroojlooy jadid
Department of Industrial and Systems Engineering, Lehigh University
Harold S. Mohler Laboratory
200 West Packer Avenue
Bethlehem, PA 18015-1582, USA
Tel: +16109741791; Cell: +14846668370
Email: afo214@lehigh.edu

[2]Department of Industrial Engineering, Sharif University of Technology, Zip code 14588/89694, Azadi Ave, Tehran, Iran,
Tel.: +982166165712; Cell: +989123873699; fax: +982166022702
Email: eshghi@Sharif.edu

Train scheduling is one of the significant issues in the railway industry in recent years since it has an important role in efficacy of railway infrastructure. In this paper, the timetabling problem of a multiple tracked railway network is discussed. More specifically, a general model is presented here in which a set of operational and safety requirements is considered. The model handles the trains overtaking in stations and considers the stations capacity. The objective function is to minimize the total travel time. Unfortunately, the problem is NP-hard and real size problems cannot be solved in an acceptable amount of time. In order to reduce the processing time, we presented some heuristic rules, which reduce the number of binary variables. These rules are based on problem's parameters such as travel time, dwell time and safety time of stations and try to remove the impracticable areas of the solution space. Furthermore, a Lagrangian Relaxation algorithm model is presented in order to find a lower-bound. Finally, comprehensive numerical experiments on the Tehran Metro case are reported. Results show the efficiency of the heuristic rules and also the Lagrangian Relaxation method in a way that for all analyzed problems the optimum value are obtained.

## 1- Introduction

Railway is a fast and economic mode of transportation; according to the Association of American Railroad's study, rail companies move more than 40 percent of the US's total freight [1] and is predicted to expand as double of the current amount on 2020. So, the railways managers have to expand the infrastructures or manage the current facilities more efficiently. The construction of the new infrastructures is very expensive and time consuming so that the utilization efficiency of the current network facilities by optimized line planning, network timetabling, crew scheduling and maintenance scheduling is very important. One of the most influential majors in this list is timetabling, which has engendered a big field of study by itself.

In a general point of view, it can be divided into three main fields: mathematical programing, simulation based optimization methods and expert systems. In the field of the mathematical programing, the aim is to create the global optimal timetable or reschedule the existing timetable. This article is focused on creating new schedule and timetables, based on mathematical programming.

[2] presented a mathematical model for single-track railway with dynamic travel times. They considered delay and operational costs as the objective function and proposed a Branch and Bound (B&B) to solve the problem. [3] presented a model based on Periodic Event Scheduling Problem (PESP) in which all of their parameters such as dwell time, headway time and trip time are dynamically considered. [4] developed a multi-objective nonlinear model to minimize fuel consumption and total passenger-time. They found the Pareto frontier and then used a distance-based method to find the solution. [5] used the Pareto solution to solve a double objective model and, a combination of the expected waiting time and the total travel time as an objective function. A beam search in a B&B algorithm is proposed to solve the MIP problem. In another study on multi-objective problems [6] proposed a Particle Swarm Optimization (PSO) for dealing with the total travel time and the variation of inter departure. [7] proposed an objective function consists of different types of waiting time and the late arrival. Furthermore, they presented a two phase algorithm, obtaining an ideal buffer time, and then a timetable is created by using an LP model. Finally, a simulation compares different timetables. [8] proposed ideal running time instead of summation of actual travel and dwell time. Their objective is based on the delay distribution of trains, the passengers count and different types of waiting time and late arrivals. They built different timetable with a LP model and a simulation evaluates them. [9] proposed a mixed integer model with a fuzzy multi objective function, to minimize energy minimization, carbon emission cost and total passenger time. The model considered an improved version of objective function compared to [4] and in the situation that all trains powered by electricity, the proposed model degenerates to [4] model.

[10] presented a complicated B&B algorithm to solve the proposed problem in [5]. They developed three methods for node selection and proposed some other complicated rules for the branching process in B&B method. Also, they designed a Lagrangian Relaxation method and another heuristic method to find a lower bound. [11] proposed a heuristic method that provides train path and timetable simultaneously on mixed single and double track networks. The method has four phases that iteratively creates and adjusts a timetable. [12] for a single line network proposed a three stage method which decomposes the model to find a solution with the maximum relative time, i.e. the ratio of travel time to minimum possible travel time, the minimum sum of the departure time, and the minimum fuel consumption. In order to solve the problem, they proposed a bicriteria algorithm to minimize the relative travel time. [13] proposed a model to minimize the total passenger trip time, considering the number of passengers as a stochastic function. They used the expected value, pessimistic and optimistic values for the number of passengers and designed a B&B algorithm to solve the model. [14] proposed a bisection method with objective function of relative travel time and used some heuristic methods to reduce the number of inactive binary variables in single track and double tracks networks. Furthermore, they proved that the solution of the algorithm is optimal. [15] presented a nonlinear model which considers alternate double–single track (ADST) lines. Their objective is to minimize the construction cost, maximum relative travel time of all trains, and the sum of total relative travel times such that the obtained departure times be close to the desired ones. They linearized the initial model and proposed some binary reduction.

By the graph approach, [16] modeled the train timetabling problem as a blocking parallel-machine job shop scheduling problem and introduced an improved Shifting Bottleneck Procedure. Then, they used an alternative graph to solve the main problem. [17] introduced a graph model for parallel rails with linked rails in sidings, capacitated buffer, acceleration and deceleration time which also does not allow unforced idle time. In order to solve the problem, they used a constructive algorithm (CA), simulated annealing (SA) and local search (LS) metaheuristics. [18] addressed the adjusting of timetables to handle perturbations and unnecessary multiple overtaking conflicts. They used the disjunctive graph to represent the problem and used LS and SA to obtain a good correction. [19] considered the timetabling problem as a job-shop problem and addressed a customized disjunctive graph to construct the timetable. In addition, they proposed some CA to create feasible solutions. [20] introduced a disjunctive graph model, considering trains and sections length, headways, and blocking conditions, non-delay scheduling policy and passing loops. They proposed a CA, based on NEH algorithm for Job shop problem, an SA, and LS methods to create and improve solutions. [21] proposed a new model for solving microscopic scale train timetabling problems. In their defined problem, stations have multi line and gate capacities and the arrival and departure time of trains are known. They proposed a new method for determining large conflict cliques in conflict graph and put them in an ILP model, relaxing strongly the related LP problem. [22] proposed a disjunctive graph, considering capacitated stations and sidings on single line networks, no-wait, and blocking properties. A hybrid algorithm is proposed to construct the feasible train timetable combined by a LS algorithm, minimizing the makespan. [23] proposed a heuristic algorithm using relaxed dynamic programming, based on the acyclic space-time graph. The algorithm starts with an ideal timetable and tries to resolve the conflicts between trains. Their model considers single corridor networks, parallel tracks in the stations, and also the presence of junctions on the network.

In the category of meta-heuristic methods, [24] proposed a Genetic Algorithm (GA) for the PESP that included a guided process to build initial population. [25] proposed a model to deal with train timetabling in single rail networks which is based on PESP. They also proposed a hybrid algorithm consisting of SA and PSO. [26] proposed a new robust periodic model which is based on PESP on a single rail network. In their model, station capacity and headway constraints are considered. A fuzzy approach is used to consider the robustness and its tradeoff with the train delay and the time interval between departures of trains from same origin. Also, a SA is used to solve the problems. In the robust problem category, [27] presented a complete survey on the robust models and their features. Some new robust research also is presented in [28]. [29] proposed a mathematical programming model, considering total energy consumption and total traversing time optimization in railway network with multiple trains and multiple links in stations. An integrated GA and simulation is used to obtain an approximate optimal strategy. [30] proposed a customized Ant Colony Optimization to minimize total weighted tardiness in single track railways. [31] proposed a travel advance strategy (TAS) method combined with GA for dealing with single track railways. The algorithm searches for a timetable with minimal delay ratio, i.e. the total delay time over the total free-run time. [32] proposed a model to obtain timetabling on one-way high speed double track networks. In order to solve the model, they proposed an improved GA and used simulation to analyze the accuracy of the algorithm. [33] also proposed a GA to provide the timetable of an urban rail transit system. The model adjusts the headway to obtain the best trade-off between the passenger travel time and energy consumption with a guaranteed transit capacity.

In a little bit different context, [34] proposed a linear formulation of cyclic timetabling problem for single track railways in which minimization of cycle length is the objective. [35] proposed a train timetabling model which deals with dynamic demand. Considering dynamic demand for different routes in different times, they proposed three binary models and a branch and cut algorithm to minimize the total waiting of passengers in stations.

We could not find any study representing a linear mathematical model to obtain timetabling in a network with parallel uni-directional tracks and limited number of station's platforms or siding's capacities. In this way we propose an approach that models overtaking decisions at stations/ sidings, as opposed to other approaches which model precedence's (i.e. sequencing) on single tracks as limited sources. Also, we show that the problem is NP-hard and the real size problems cannot be solved in a reasonable amount of time. In order to obtain beneficial solutions in a reasonable time, we provide some new upper bound and lower bound rules. By concentrating on these two topics, the structure of this article is as follows. In section 2 we provide the mathematical model. Then in section 3 we propose the upper bound rules and in section 4 propose the lower bound rule and it is followed by numerical experiment's result in a real world problem.

## 2- Problem Definition

In this section, we define the problem and its assumptions, give the notation and introduce our model.

This study considers a general situation of a double-lines networks. In all lines and corridors in the network, there are two parallel tracks and each track facilitates uni-directional flow and not bi-directional, i.e. there is no opposite direction train, which is a common situation in metro, subway and high speed networks; and is of interest as [36, 37] covered some earlier research. Each station or siding can have one or more capacity, which some of them may have no platform for passengers. Also, in sidings there is no linked rail between opposite directions and each direction has its separate sidings. In addition, overtaking is allowed and trains can only overtake each other in the siding (passing loops) or in the stations with more than one platform, i.e. train routes in stations are not fixed. Moreover, the main line is the line that connects the stations or siding with each other. A general view of siding and main line are showed in Figure 13.

[Figure 1 should be here]

There may be different types of trains with different speeds and dwell time. The trains may be freight or passenger specific or any mix of them. It is assumed that for each train in each station of the network, there is a minimum and maximum dwell time and also travel time between two consecutive stations or sidings is known. Routes and the dispatching sequence of trains between stations and sidings are given. In addition, the travel time between two stations or sidings is the difference between the departure times of the last station until the arrival time of the current station. Considering the type of trains and their intrinsic nature in metro and high speed networks, there is no significant time for acceleration and deceleration time, but we have considered the minimum travel time between two stations to include these times. Note that, the considered parameters, e.g. travel time, dispatching sequence, and dwelling time are known in any railway network and they are the most common features of the MIP-models for the timetabling problem.

A mathematical model is presented as follows. In this model output is the timetable of a railway network. The most important assumptions of our model are as below.

- All parameters of the model are deterministic.
- Each corridor of the network has one line in each direction.
- The station capacity which is the number of platforms in the stations can be any positive integer number. There is a same situation about sidings.
- Acceleration and deceleration time are considered in the travel time.

- Number of trains and their routes are given.

2-1  Notation

Throughout the paper, we reserve $l$ to denote line index in a network, $L$ is the list of all lines, $t$ defines train index , $T_l$ is the list of trains in line $l$, $s$ is the station index in each line and $S_l$ defines the total number of station(s) in line $l$. Also, $q$ is the next station in line $l$, .i.e. $q = s + 1$. Finally, $m_l$ denotes the last station in line $l$. The major parameters of our problem are as follows:

$SF_{l,s}$: Minimum headway (Safety time) between two adjacent trains at station $s$ at line $l$.

$SC_s$: Capacity of station $s$, which is the number of platform tracks and capacity of stops.

$\underline{D}_{l,s}^t$: Minimum dwell time of train $t$ in station $s$ of line $l$.

$\overline{D}_{l,s}^t$: Maximum dwell time of train $t$ in station $s$ of line $l$

$\underline{T}_{l,s,q}^t$: Minimum travel time of train $t$ between station $s$ and $q$ $(s+1)$ of line $l$

$\overline{T}_{l,s,q}^t$: Maximum travel time of train $t$ between station $s$ and $q$ $(s+1)$ of line $l$

$r_{l,t,s}$: Earliest start time. Minimum start time of train $t$ in station $s$ of line $l$

$M$ : shows a large positive number.

Also, we use these two notations throughout the paper:

$\psi_{t,t',q}$: The time interval between the departure time of train $t$ from station $q$ and the arrival time of train $t'$ to this station

$\vartheta_{t,t',s}$: The time interval between the departure time of two trains from station $s$

For safety reasons trains are not permitted to close to each other in case of collision. Between departure and arrival of trains in each station or siding minimum headway must be satisfied. Other situations are, departure of two trains from a station or siding and arrival of two trains in a station or siding. Stations' headway is a predefined parameter that is inherently related to geographical and technical characteristics of each station. In order to identify the decision variables, first, we define the situations that every two trains may encounter in a station or siding. These situations are the events that may bring some changes to the sequence of any two arbitrary trains. With regard to the situations of the problem, for each two trains, four different situations can occur. Event (1) defines the situation that train $t$ overtakes train $t'$ in station $s$. In event (2) train $t'$ overtakes train $t$ in station $s$. In event (3) train $t'$ arrives and departs station $s$ after train $t$. Finally, event (4) demonstrates the situation that train $t'$ arrives and departs station $s$ before train $t$. In other words, in this event train $t$ has overtaken train $t'$ before station $s$ or is scheduled before train $t'$ from the first station. Figure 14 shows and clarifies the mentioned situations for the trains $t$ and $t'$.

[Figure 2 should be here]

By considering these events and situations, the decision variables are as below that cover any situation in a network:

$c_{l,t,s}$ = Departure time of $t^{\text{th}}$ train from station $s$ of line $l$
$s_{l,t,s}$ = Arrival time of $t^{\text{th}}$ train in station $s$ of line $l$

$x_{t,t',s}$ = 1 if train $t'$ overtakes train $t$ in station s as event 1 and 2, otherwise 0

$x'_{t,t',s}= 1$ if train $t'$ scheduled after train $t$ in station $s$, otherwise 0. This variable represents the situations of event 3 and 4.

$y_{t,t',s} = 1$ if train $t'$ and train $t$ use simultaneously station s, otherwise 0

As the definition of decision variables shows, binary variables will model overtaking decisions in a double-line uni-directional railway network. By this approach, our decision variables define the sequence of trains in the network and determine when and in which station or siding the trains will overtake each other. Somehow, this is a different approach compared to the current models that define precedence and sequencing on single tracks. With regard to these variables, constraints of problem are defined as follows.

### 2-2 Objective Function

Our model's aim is to minimize the summation of the departure time at the last station for all trains, which is the minimum total travel time as is shown in Equation (1).

$$\min z = \sum_{l=1}^{L} \sum_{t=1}^{T} c_{ltm_l}; \tag{1}$$

This objective function tries to increase the efficacy of the network's infrastructure and also minimize the travel time. Therefore, it satisfies the criteria of two main groups in the railway networks, i.e. passengers and railway owners. It is also a commonly used objective function in research of [3, 5, 10, 12, 36, 38, 39].

### 2-3 Travel Time Constraints

The constraints below ensure that the running time of a train does not violate a given limit of time, which also defines the speed limitation in the line. Furthermore, they do not allow the running time to drop below a specific amount of time defined by the desire of passengers. For all $t \in T_l; \ \forall \ l \in L; \forall \ s \in S_l$ the constraints are as follow:

$$s_{l,t,s+1} \geq c_{lts} + \underline{T}^t_{l,s,q} \tag{2}$$

$$s_{l,t,s+1} \leq c_{lts} + \overline{T}^t_{l,s,q} \tag{3}$$

### 2-4 Dwell Time Constraints

The following constraints ensure that the dwell time of a train does not violate a given limit of time. The upper and lower bounds must be gathered as real values for each station and train. Moreover, for trains that have not been scheduled to stop at a given station, both bounds are set to zero for that station. For all $t \in T_l; \ \forall \ l \in L; \forall \ s \in S_l$ the constraints are as below:

$$c_{lts} \geq s_{lts} + \underline{D}^t_{l,s} \tag{4}$$

$$c_{lts} \leq s_{lts} + \overline{D}^t_{l,s} \tag{5}$$

### 2-5 Overtaking Constraints

The constraints below aim to define the sequence of trains in each corridor and station and also determine the arrival and departure time of trains. Using the defined binary variables, related constraints to each event have been defined. In event (3) and (4), constraints (6) and (7) ensure that the arrival and departure time of train $t'$ be greater than train $t$ if train $t'$ scheduled after train $t$ and vice versa. Also, in this situation constraint (8) does not allow trains overtake each other except in the stations. Constraint (9) and (10) ensure that the arrival time of train $t'$ be greater than train $t$ and its departure be lower than train $t$, if train $t'$ scheduled for overtaking of train $t$. They act the same as description of events (1) and (2). Constraint (11) acts similarly to the constraint (8) and does not allow in events (1) and (2) trains overtake each other except in the stations. Constraint (12) ensures that in each station and for each pair of trains $t$ and $t'$, exactly one event can occur. In order to simplify the notation, we used $\phi'_{l,s}(t, t')$ to represent $\text{SF}_{l,s} - \left(1 - x'_{t,t',s}\right)M$. For all $l \in L;\ \forall\, t, t' \in T_l;\ \forall\, s \in S_l$ the constraints are as follow.

$$s_{lt's} \geq s_{lts} + \phi'_{l,s}(t, t') \quad ;\ \forall\, t \neq t' \tag{6}$$

$$c_{lt's} \geq c_{lts} + \phi'_{l,s}(t, t') \quad ;\ \forall\, t \neq t' \tag{7}$$

$$s_{lt's} \geq s_{lts} + \phi'_{l,s}(t, t') \quad ;\ s > 1 \quad ;\ \forall\, t \neq t' \tag{8}$$

$$s_{lt's} \geq s_{lts} + \phi_{l,s}(t, t') \quad ;\ \forall\, t \neq t' \tag{9}$$

$$c_{lts} \geq c_{lt's} + \phi_{l,s}(t, t') \quad ;\ \forall\, t \neq t' \tag{10}$$

$$s_{lts} \geq s_{lt's} + \phi_{l,s}(t, t') \quad ;\ s > 1 \quad ;\ \forall\, t \neq t' \tag{11}$$

$$x'_{t',t,s} + x'_{t,t',s} + x_{t,t',s} + x_{t',t,s} = 1 ;\ \forall\, t \neq t' \tag{12}$$

### 2-6  Station Capacity Constraints

This set of constraints is proposed to guarantee the maximum capacity constraint of each station or siding. Constraint (13) assures that in a station with one capacity, only one train stops in station at any time (event type (3) and (4)). Thus the arrival time of one train must be greater than another one which results in the condition that more than one train cannot be stopped at one moment in the station.

$$s_{lts} \geq c_{lt's} + \phi'_{l,s}(t', t)$$
$$s_{lt's} \geq c_{lts} + \phi'_{l,s}(t, t') \qquad ;\ \text{SC}_s = 1 \quad ;\ \forall\, t > t' \tag{13}$$

$$\left(x_{t_1,t_2,s} + x_{t_2,t_1,s}\right) + \left(x_{t_1,t_3,s} + x_{t_3,t_1,s}\right) + \cdots + \left(x_{t_1,t_{\text{SC}_s+1},s} + x_{t_{\text{SC}_s+1},t_1,s}\right) + \cdots +$$
$$\left(x_{t_{\text{SC}_s},t_{\text{SC}_s+1},s} + x_{t_{\text{SC}_s+1},t_{\text{SC}_s},s}\right) \leq \binom{\text{SC}_s+1}{2} - 1 ; \begin{cases} \forall\, l \in L \\ \forall\, t_1, t_2 \ldots t_{\text{SC}_s+1} \in T_l \\ \forall\, s \in S_l\ and\ \text{SC}_s > 1 \end{cases} \tag{14}$$

Constraint (14) considers the capacity constraint of station $s$ with capacity of $SCs$, for the overtaking events (event type (1) and (2)). In stations that their capacity is one, constraint (13) does not allow any overtaking, i.e. the constraint is $\left(x_{t,t',s} + x_{t',t,s}\right) \leq 0$ for all of $t \neq t' \in T_l$ and $\forall\, s \in S_l$ in line $l$. In stations with one more capacity, the constraint considers all pair combinations of three trains (because $\binom{SC_s+1}{2} = 3$) as $\left(x_{t,t',s} + x_{t',t,s}\right) + \left(x_{t,t'',s} + x_{t'',t,s}\right) + \left(x_{t',t'',s} + x_{t'',t',s}\right) \leq 2$ for all $t \neq t' \neq t'' \in T_l$ and $\forall\, s \in S_l$ in line $l$. This constraint ensures at most one train can overtake from another train at each time, i.e. at each time only two trains can stop in the station $s$. To further illustrate the issue, consider the situation that in station $s$, which has two platforms, train $t$ is scheduled to overtake train $t'$, i.e.

$x_{t',t,s} = 1$, and $(x_{t,t',s} + x_{t',t,s}) = 1$, and train $t''$ is entering the station $s$. Three situations may occur:

1- Assume train $t''$ is scheduled to depart the station $s$ after train $t'$. In this situation there is no additional overtaking and train $t''$ may incur some delay.

2- Assume that without regarding the station capacity, train $t''$ is scheduled to overtake train $t$ in the station $s$.

If train $t''$ overtakes train $t$, i.e. $x_{t,t'',s} = 1$ and $(x_{t,t'',s} + x_{t,t,s}) = 1$, train $t''$ arrive later than trains $t$ and $t'$ and must depart before them. In another word, train $t''$ also has to overtake train $t'$ that means $x_{t',t'',s} = 1$ and $(x_{t',t',s} + x_{t',t'',s}) = 1$. So, $(x_{t,t',s} + x_{t',t,s}) + (x_{t,t'',s} + x_{t'',t,s}) + (x_{t',t'',s} + x_{t'',t',s}) = 3$ and considering the RHS of the constraint, the situation is not allowed and the assumption is not true. So, constraint (14) does not allow exceeding the capacity of the station.

3- Assume train $t''$ is scheduled to overtake train $t'$ in the station $s$ and train $t$ depart the station before arrival of train $t''$.

If train $t''$ and train $t$ overtake train $t'$ in station $s$, i.e. $x_{t',t'',s} = 1$ and $(x_{t'',t',s} + x_{t',t'',s}) = 1$ and $x_{t',t,s} = 1$ and $(x_{t,t',s} + x_{t',t,s}) = 1$ the constraint satisfies.

With a similar approach for station with capacity of $SC_s$, there are $\frac{(SC_s+1)SC_s}{2}$ pairs of variables and on the right hand side, the number of pairs minus one is replaced. This set does not allow violation in stations' capacity with a larger capacity. Furthermore, constraints (15), (16) and (17) consider the station capacity when the two trains are scheduled with events type (3) or (4).

$$s_{lt's} \geq c_{lts} + \text{SF}_{l,s} - (1 - x'_{t,t',s})M - My_{t,t',s} \quad ; \quad SC_s \geq 2 \ ; \ \forall t > t' \tag{15}$$

$$\sum_{t',t\neq t'}^{T}(y_{t,t',s} + y_{t',t,s} + x_{t,t',s} + x_{t',t,s}) \leq SC_s - 1 \ ; \ SC_s \geq 2 \tag{16}$$

$$y_{t,t',s} \leq x'_{t,t',s} \quad ; \quad SC_s \geq 2 \ ; \ \forall t > t' \tag{17}$$

As mentioned $y_{t,t',s}$ is the auxiliary variable to consider the station capacity in the mentioned events. Constraints (15) and (17) insure that the auxiliary variables can take value when event types (3) and (4) occur. For instance in event type 3 that train $t$ scheduled before train $t'$, i.e. $x'_{t,t',s} = 1$, $y_{t,t',s}$ can take value one. Regarding Constraint (15), if $y_{t,t',s}$ takes 1, $s_{lt's}$ will be greater than $-M$, this means that the trains can stop simultaneously in station $s$. The same situations can be assumed for event type (4). Also, Constraint (16) insures that the number of trains that can stop in each station will not be more than its capacity.

These constraints with continues time variables handle the station capacity in each of the four events, and it is a new formulation in the literature. The current papers mostly consider time slotting approach, e.g. [41], which consider a time horizon for the network's operation time and discretize it into some timestamps and propose some constraints to consider the station capacity. A similar slot based scheduling is used by [39] in an integrated bi-level station layout design and scheduling model.

### 2-7  Departure Time Constraint

It is necessary that some trains may have a predefined departure time for some stations. It may occur in the first station of the corridor. For all $l \in L; \ \forall t \in T_l; \forall s \in S_l$ the constraint is as equation (18).

$$c_{lts} \geq r_{lts} \tag{18}$$

Compared to the current models of train timetabling problem in the literature, our model has some advantageous points. The first feature of our model is to define event base decision variables and constraints. These variables make it possible to use them in decompositions or relaxation algorithms so that each sub-problem consists of a decision making problem related to each event. Constraint (12) which integrates the decision variables allows the occurrence of only one event and creates a suitable situation for decomposition. In this situation a difficult problem can be decomposed to some simpler sub-problems and increases the possibility of solving real size problems.

On the other hand, as another benefit of the model, planning of the station capacity is simultaneously possible in our timetabling model. Constraints (14) and (16) make it possible to consider the capacity of stations and sidings. Most of models in the literature decompose the problem to two separate problems. One to define and solve a general timetable and another to deal with the complex station capacity planning problem [2, 5, 10, 12, 14, 37, 42, 43, 44]. The station capacity planning problem needs special attentions, so some research and algorithms in recent years focused on this field of study [21, 45, 46, 47, 48, 49, 50, 51]. By this point of view, the proposed model integrated the two models and provides an integrated answer. Moreover, the new definition of decision variables classifies the solution space and gives possibility to define and remove the areas of solution space that are not rationally as part of optimal solution area. A comprehensive investigation of this topic will follow.

## 2-8  Complexity Analysis

Although the problem is generally considered as a NP-hard problem, almost every paper considers a different version of the problem. For this reason, we give a NP-hardness proof for the specific problem we consider. Our proof will show that our problem is a reduction of $FFs/ \ r_J \ / \sum C_j$. In order to show the complexity of the problem, we used a reduction from a flexible flowshop model. Our proposed model is a complicated version of a flexible flow-shop model with capacitated buffer and some limitations for the job assignment to machines. According to the definition of flowshop problem by [52], there are $m$ machines in series and each job passes each machine in a same route. If the number of identical machines in at least one stage be greater than 1, the problem is classified as a Flexible Flow Shop problem ($FF_s$). Moreover, when machine $M$ in a specific stage is not capable for processing all jobs, the $M_j$ defines this machine eligibility constraint. Finally, if there is a limited buffer between two machines, in a way that when the buffer is full the machine is not allowed to release a new job, the *block* constraint obligates the completed job to remain in the previous machine.

To further illustrate the issue, consider the trains as jobs and the stations as machines which each job (train) must pass through machines (stations); as [10, 20] have proposed. As the definition of flowshop problem by [52] this problem can be classified as flow-shop problem with *block* constraint. Moreover, the number of platforms in each station or siding defines the capacity of the station or siding that may be greater than one, so the problem classifies as a flexible flowshop problem or multi-stage parallel machine problem [20]. On the other hand, some express trains are scheduled to not stop in some stations and the stations also have some specific platforms for non-stopping trains. This situation defines the machine eligibility constraint. Furthermore, there are intrinsically some limitations about the start time of each train. Thus, the problem can be shown as $FF_S/ \ M_J; \ r_J;$ block$| \ \gamma$ according to the notation of [52]. Since the problem of $FF_S//\sum C_j$ is NP-hard [53], our problem is NP-hard. Therefore, it is not possible to solve large scale problems in a poly-nominal time and some policies should be conducted to reduce complexity of the problem.

### 3- Upper Bound Rules

In this section, we define some rules to efficiently limit the number of binary variables and also create efficient upper bounds. Binary variables, which increase the complexity of the problem, are used to identify the sequence and overtaking details of trains. Therefore, decreasing the number of binary variables decreases the number of decisions and make the problem easier to solve. There exists some research in order to limit the number of unnecessary over-takings. [18] introduced some algorithms to identify and correct some conflicts such as multiple overtaking and compound moves. In their article they restrict trains from multiple overtaking and there is no point about restricting general situations of overtaking. In the following we propose some general rules to restrict undesirable overtaking among the trains in a network. We use trains' parameters in order to define usefulness or undesirability of overtaking among every two trains.

In each railway network, some trains have same characteristics such as travel and dwell time. It is clear that the study of this group of trains can lead to rules that may relax the binary decision variables among the similar trains. For example, in German railway network there are eight types of train such as Intercity Express (ICE), RegionalBahn (RB) and S-Bahn (S) and most trains in each group have same parameters in the network. Consider two trains $A$ and $B$ in a group that have exactly same speed and dwell time parameters in all network stations. The idea is that overtaking of each of these two trains from each other is not beneficial. Consider the timetable that train $A$ is scheduled to overtake train $B$ in station $s$. Because the dwell time and travel time of two trains are identical, train $B$ has to delay at the station and as a result the objective function increases. By this idea, we try to expand the simple rule to the whole of the trains and considering the general situations we will obtain the rules that fix unnecessary overtaking binary variables. In a word, we want to obtain the situations among every two trains that occurrence of events (1) and (2) may result in a delay in the timetable. We consider two random trains and the objective function of two timetables in order to find the situations and related rules that overtaking is not beneficial. We compare the objective functions of the timetable that trains go ahead without overtaking with the timetable that one train overtakes another one. In these investigations, we assume that there is no deliberate delay. The results of the investigations can be used in reduction of binary variables and the CPU time.

Now, assume two trains $t'$ and $t$, which are now in stations $s$ and $q$ ($= s+1$) consecutively as shown in Figure 15, which station $q$ has two platforms. In the Figure 15, the two horizontal lines show the scheduling of the related stations. Points $k$ and $w$ define the arrival and departure time of train $t$ in station $q$ respectively. The time interval between $R$ and $K$ represents travel time of train $t$ between two consecutive stations $s$ and $q$. Also, the time interval between $K$ and $W$ defines the dwell time of train $t$ in station $q$.

[Figure 3 should be here]

At first, we define two time intervals in order to characterize different situations. We define $\psi_{t,t',q}$ as the time interval between the departure time of train $t$ from station $q$ and the arrival time of train $t'$ to this station, i.e. $\psi_{t,t',q} = S_{tl'q} - C_{tlq}$. Also, $\vartheta_{t,t',s}$ is the time interval between the departure time of two trains from station $s$, i.e. $\vartheta_{t,t',s} = C_{lt's} - C_{lts}$. If $\psi_{t,t',q}$ be greater than zero, the overtaking of train $t'$ from train $t$ in station $q$ is not recommended, because overtaking will result in a delay for train $t$. So, the timetable without overtaking is selected for the station and trains. Moreover, without loss of generality, we can assume that $\vartheta_{t,t';s}$ is greater than zero. The mentioned condition can be written as follows:

If $0 < \psi_{t,t';q}$ and $0 < \vartheta_{t,t';s}$ then $x_{t,t';s} = 0$.

In some situations this rule does not give the optimal solution, although it fixes lots of binary variables. To investigate further, we analyze the situation that the mentioned condition is not satisfied, i.e. two trains simultaneously are in the station and one train can overtake another one, that means $\psi_{t,t';q}$ be lower than zero. This condition meets when the time interval between the arrival and departure of two trains in station $q$ is lower than the required station safety time. The first derived condition is:

$$\psi_{tt'q} \leq 0 \; \equiv \; S_{lt'q} \leq C_{ltq} + SF_{l,q} \tag{19}$$

$$C_{lt's} + \underline{T}_{l,s,q}^{t'} \; \leq \; S_{lt'q} \leq (C_{lts} + \underline{T}_{l,s,q}^{t} + \underline{D}_{l,q}^{t}) + SF_{l,q} \tag{20}$$

$$\equiv \; C_{lt's} - C_{lts} \; \leq \underline{T}_{l,s,q}^{t} - \underline{T}_{l,s,q}^{t'} + \underline{D}_{l,q}^{t} + SF_{l,q} \tag{21}$$

In the equation (20), it is obvious that, the departure time of train $t$ from station $q$ is equal to its departure time from station $s$ plus its travel time between two stations and its dwell time in the station $q$, i.e. $C_{ltq} = C_{lts} + \underline{T}_{l,s,q}^{t} + \underline{D}_{l,q}^{t}$. Also, as another fact, the arrival time of train $t'$ to station $q$ is greater than its departure time from station $s$ and the travel time between two stations, i.e. $S_{lt'q} \geq C_{lt's} + \underline{T}_{l,s,q}^{t'}$. On the other hand, it is obvious that $\vartheta_{t,t';s}$ is greater than dwell time of train $t'$ in station $s$ plus safety time of this station, i.e. $\underline{D}_{l,s}^{t'} + SF_{l,s} \leq C_{lt's} - C_{lts}$. Thus, we can write:

$$\underline{D}_{l,s}^{t'} + SF_{l,s} \leq C_{t's} - C_{ts} \leq \underline{T}_{l,s,q}^{t} - \underline{T}_{l,s,q}^{t'} + \underline{D}_{l,q}^{t} + SF_{l,q} \tag{22}$$

The derived condition in equation (22) describes the situation that two trains are simultaneously in the station $q$ so that overtaking may be beneficial and we will investigate it later on. On the other hand, if the condition in equation (22) does not meet (like the situation that is shown in Figure 15), overtaking of train $t'$ from train $t$ in station $q$ results at least $\underline{D}_{l,s}^{t'} + SF_{l,s} - (\underline{T}_{l,s,q}^{t} - \underline{T}_{l,s,q}^{t'} + \underline{D}_{l,q}^{t} + SF_{l,q})$ obligatory delay for train $t$. This situation causes the increasing of objective function and a grasp vision constraint (23) can be added to the problem under the below conditions:

$$\forall \, l \, \in \, L \; ; \; \forall \, s,q \; \in S_l ; \forall \, t \neq t' \in \, T_l$$

$$\underline{D}_{l,s}^{t'} \; + \; SF_{l,s} > \underline{T}_{l,s,q}^{t} - \underline{T}_{l,s,q}^{t'} + \underline{D}_{l,q}^{t} + SF_{l,q}$$

and the constraint is as follows that forbids overtaking of two trains in station $s$.

$$x_{t,t',s} \leq 0; \tag{23}$$

In a specific condition the added constraint brings the optimal solution as follows.

**Proposition 1**: Suppose that there is a railway network with three stations and no opposite direction trains. Furthermore, there are two consecutive trains, $t$ and $t'$, that at the first station train $t$ precedes train $t'$. Train $t$ will finish its path first, if $\vartheta_{t,t';1}$ be greater than $\psi_{t,t';2}$, i.e. the time interval between the departure time of two trains in first station ($\vartheta_{t,t';1}$) be greater than the time interval between the arrival time of train $t'$ to the second station and the departure time of train $t$ from second station ($\psi_{t,t';2}$).

**Proof**: As described in the proposition, there are three stations and the only station where overtaking is possible is the second station. Moreover, as described in the proposition, the condition in equation (22) is not satisfied. If train $t'$ overtakes train $t$ in the second station, as shown in equation (22), the objective function at least increases by $\underline{D}_{l,s}^{t'} + \mathrm{SF}_{l,s} - (\underline{T}_{l,s,q}^{t} - \underline{T}_{l,s,q}^{t'} + \underline{D}_{l,q}^{t} + \mathrm{SF}_{l,q})$. So the optimal route is the current route that train $t$ ends its route first. $\blacksquare$

On the other hand, we consider the situation that the condition in equation (22) meets, i.e. $\psi_{t,t',q} < 0$ and $0 < \vartheta_{t',s}$ and in a word, $\underline{D}_{l,s}^{t'} + \mathrm{SF}_{l,s} \leq \underline{T}_{l,s,q}^{t} - \underline{T}_{l,s,q}^{t'} + \underline{D}_{l,q}^{t} + \mathrm{SF}_{l,q}$. In this situation train $t'$ arrives in station $q$ before train $t$ departs the station. In order to simplify the investigation and also make the equations more clear, we try to classify different circumstances and then check their conditions. As the first classification factor, we consider the relation between $\psi_{t,t',q}$ and $\mathrm{SF}_{l,q}$. In the situation that the equation (22) is valid, we consider the condition that $|\psi_{t,t',q}|$ be lower than the station's safety time, i.e. $\psi_{t,t',q} < \mathrm{SF}_{l,q}$, so train $t'$ will be delayed until the safety time meets. With regard to this condition two different situations in which the safety time can/ cannot satisfy will be discussed. In order to obtain the relation between $\psi_{t,t',q}$ and $\mathrm{SF}_{l,q}$ more clearly, we use equation $\psi_{t,t',q} = C_{lt's} + \underline{T}_{l,s,q}^{t'} - C_{lts} - \underline{T}_{l,s,q}^{t}$ which describes the time interval between arrival time of two trains to station $q$ and according to inequality $\underline{D}_{l,s}^{t'} + \mathrm{SF}_{l,s} \leq C_{t's} - C_{ts}$ is written as $\psi_{t,t',q} = \mathrm{SF}_{l,s} + \underline{D}_{l,s}^{t'} + \underline{T}_{l,s,q}^{t'} - \underline{T}_{l,s,q}^{t}$. So, the first classification can be written as $\mathrm{SF}_{l,s} + \underline{D}_{l,s}^{t'} + \underline{T}_{l,s,q}^{t'} - \underline{T}_{l,s,q}^{t} \leq \mathrm{SF}_{l,q}$ or $\mathrm{SF}_{l,s} + \underline{D}_{l,s}^{t'} + \underline{T}_{l,s,q}^{t'} - \underline{T}_{l,s,q}^{t} > \mathrm{SF}_{l,q}$.

Moreover, as the second classification factor, we consider the situation that train $t'$ overtakes train $t$ in station $q$ which causes scheduled delays for train $t$. In another situation that there is not any scheduled delay for train $t$, i.e. the time interval that train $t$ stops in the station $q$ is greater than the required safety time for train $t'$ to enter the station $q$ ($\mathrm{SF}_{l,q}$) plus its dwell time in that station ($\underline{D}_{l,q}^{t'}$) plus another required safety time ($\mathrm{SF}_{l,q}$) after train $t'$ departs station $q$. So, $\mathrm{SF}_{l,q} + \underline{D}_{l,q}^{t'} + \mathrm{SF}_{l,q} \geq \underline{D}_{l,q}^{t}$ and $\mathrm{SF}_{l,q} + \underline{D}_{l,q}^{t'} + \mathrm{SF}_{l,q} < \underline{D}_{l,q}^{t}$ are the two classifications. According to these two classifications, four different conditions are obtained, their characteristics and assumptions are summarized in (24). Also, one example from each of them is shown in Figure 16. For all $l \in L; \ \forall \ t \neq t' \ \in T_l; \ \forall \ s,q \ \in S_l$ the conditions are as:

$$
\begin{aligned}
&\underline{D}_{l,s}^{t'} + \mathrm{SF}_{l,s} \leq \underline{T}_{l,s,q}^{t} - \underline{T}_{l,s,q}^{t'} + \underline{D}_{l,q}^{t} \\
&+ \mathrm{SF}_{l,q}
\begin{cases}
\mathrm{SF}_{l,s} + \underline{D}_{l,s}^{t'} + \underline{T}_{l,s,q}^{t'} - \underline{T}_{l,s,q}^{t} \leq \mathrm{SF}_{l,q}
\begin{cases}
\underline{D}_{l,q}^{t'} + 2\mathrm{SF}_{l,q} \geq \underline{D}_{l,q}^{t} \ (1) \\
\underline{D}_{l,q}^{t'} + 2\mathrm{SF}_{l,q} < \underline{D}_{l,q}^{t} \ (2)
\end{cases} \\
\mathrm{SF}_{l,s} + \underline{D}_{l,s}^{t'} + \underline{T}_{l,s,q}^{t'} - \underline{T}_{l,s,q}^{t} > \mathrm{SF}_{l,q}
\begin{cases}
\underline{D}_{l,q}^{t'} + 2\mathrm{SF}_{l,q} \geq \underline{D}_{l,q}^{t} \ (3) \\
\underline{D}_{l,q}^{t'} + 2\mathrm{SF}_{l,q} < \underline{D}_{l,q}^{t} \ (4)
\end{cases}
\end{cases}
\end{aligned}
\tag{24}
$$

[Figure 4 should be here]

Similar to proposition 1, for each condition we try to find the situations that overtaking is not beneficial. So, we should compare the objective functions of the two timetables that in the first, train $t'$ overtakes train $t$ and in the second, the situation that no overtaking occurs. In these comparisons, the aim is to find the relation of the parameters, in a way that overtaking by a grasp view results in increasing of the objective function.

Therefore, we obtain the objective function of the timetable that train $t$ arrives and departs stations first. We will use this amount in each of the four conditions in order to make comparisons. Then, we calculate the objective function in each of the mentioned classified conditions in a way that train $t'$ overtakes train $t$ in station $q$. These objective functions will be used in the mentioned comparison in order to obtain some heuristic rules. In the following three equations, we calculate the objective function of the no overtaking situation.

$$C_{ltq} = C_{lts} + \underline{T}_{l,s,q}^{t} + \underline{D}_{l,q}^{t} \tag{25}$$

$$C_{lt'q} = C_{lts} + (C_{lt's} - C_{lts}) + \underline{T}_{l,s,q}^{t'} + (\underline{T}_{l,s,q}^{t} - \underline{T}_{l,s,q}^{t'} + \underline{D}_{l,q}^{t} + SF_{l,q} - (C_{lt's} - C_{lts})) + \underline{D}_{l,q}^{t'} \tag{26}$$

$$C_{ltq} + C_{lt'q} = 2C_{lts} + 2\underline{T}_{l,s,q}^{t} + 2\underline{D}_{l,q}^{t} + SF_{l,q} + \underline{D}_{l,q}^{t'} \tag{27}$$

Equation (25) and (26) respectively show the departure time of trains $t$ and $t'$ from station $q$ and equation (27) is the summation of the two departure times. In equation (25), $C_{ltq} = C_{lts} + \underline{T}_{l,s,q}^{t} + \underline{D}_{l,q}^{t}$ defines the departure time of train $t$ from station $q$ without any delay. In equation (26), the departure time of train $t'$ consists of the departure time from station $s$, the travel time between two stations, the amount of time that train $t'$ must be delayed until train $t$ departs station $q$, i.e. $\underline{T}_{l,s,q}^{t} - \underline{T}_{l,s,q}^{t'} + \underline{D}_{l,q}^{t} + SF_{l,q} - (C_{lt's} - C_{lts})$ and its dwell time at station $q$. The scheduled delay amount was obtained in equation (22). As shown in equation (27), the result is $C_{ltq} + C_{lt'q} = 2C_{lts} + 2\underline{T}_{l,s,q}^{t} + 2\underline{D}_{l,q}^{t} + SF_{l,q} + \underline{D}_{l,q}^{t'}$ and will be used latter for the comparisons.

On the other hand, we have to calculate the objective function of the timetable that overtaking occurs for each of the mentioned conditions. Regarding the first set of assumptions in equation (24), the departure time of the two trains calculated in the situation that train $t'$ overtakes train $t$ in station $q$. The departure times of trains $t$, $t'$ are shown in equations (28) and (29) respectively and the objective function is as equation (30). In order to clarify the origin of the equations, the elements of each departure time is described as below:

Departure time of train $t$ = its arrival time to station $q$ + safety time of the station $q$ + dwell time of train $t'$ in the station $q$ + safety time of station $q$

Departure time of train $t'$ = its departure time from station $s$ + travel time among two station $s$ and $q$ + obligatory delay to insure the safety time of station $q$ + dwell time of train $t'$ in station $q$

$$C_{ltq} = (C_{lts} + \underline{T}_{l,s,q}^{t}) + SF_{l,q} + \underline{D}_{l,q}^{t'} + SF_{l,q} \tag{28}$$

$$C_{lt'q} = (C_{lts} + (C_{lt's} - C_{lts})) + (\underline{T}_{l,s,q}^{t'} + SF_{l,q} - (C_{lt's} + \underline{T}_{l,s,q}^{t} - C_{lts} - \underline{T}_{l,s,q}^{t})) + \underline{D}_{l,q}^{t'} \tag{29}$$

$$C_{ltq} + C_{lt'q} = 2C_{lts} + 2\underline{T}_{l,s,q}^{t} + 3SF_{l,q} + 2\underline{D}_{l,q}^{t'} \tag{30}$$

By comparing equation (27) and (30), equation (32) shows the relation of parameters in which overtaking is not recommended and the objective function of the timetable with overtaking is greater than the one without overtaking. Equation (31) shows the detail of the comparison.

$$2C_{lts} + 2\underline{T}_{l,s,q}^{t} + 2\underline{D}_{l,q}^{t} + SF_{l,q} + \underline{D}_{l,q}^{t'} \leq 2C_{lts} + 2\underline{T}_{l,s,q}^{t} + 3SF_{l,q} + 2\underline{D}_{l,q}^{t'} \tag{31}$$

$$2\underline{D}_{l,q}^{t} \leq \underline{D}_{l,q}^{t'} + 2SF_{l,q} \tag{32}$$

Thus, in the circumstance that the set of assumptions of the first condition is true and the inequality (32) is valid, we can propose that overtaking is not a good decision and we recommend adding constraint (33) to the problem under the below conditions:

$\forall\ t \neq t' \in T_l;\ \forall\ l \in L;\ \forall\ s, q \in S_l;$

$\underline{D}_{l,s}^{t'} + \mathrm{SF}_{l,s} \leq \underline{T}_{l,s,q}^{t} - \underline{T}_{l,s,q}^{t'} + \underline{D}_{l,q}^{t} + \mathrm{SF}_{l,q}$

$\mathrm{SF}_{l,s} + \underline{D}_{l,s}^{t'} + \underline{T}_{l,s,q}^{t'} - \underline{T}_{l,s,q}^{t} \leq \mathrm{SF}_{l,q};\ \ 2\mathrm{SF}_{l,q} + \underline{D}_{l,q}^{t'} \geq \underline{D}_{l,q}^{t};\ \ 2\underline{D}_{l,q}^{t} \leq 2\mathrm{SF}_{l,q} + \underline{D}_{l,q}^{t'}$

and the constraint is as follows that forbids overtaking of two trains in station $s$.

$$x_{t,t',s} \leq 0\ ; \tag{33}$$

By a similar approach, the derived equation in condition (2) is $\underline{D}_{l,q}^{t'} \leq 0$, in condition (3) is $0 \leq \underline{D}_{l,q}^{t'} + \mathrm{SF}_{l,q}$ and in condition (4) is $\underline{D}_{l,q}^{t} \leq \underline{D}_{l,q}^{t'}$. Related constraints are described in equations (34), (35) and (36) respectively. For all $l \in L;\ \forall\ t \neq t' \in T_l;\ \forall\ s, q \in S_l$ and in the condition that the parameter of problem satisfies $\underline{D}_{l,s}^{t'} + \mathrm{SF}_{l,s} \leq \underline{T}_{l,s,q}^{t} - \underline{T}_{l,s,q}^{t'} + \underline{D}_{l,q}^{t} + \mathrm{SF}_{l,q}$, the constraints are as follow.

$$x_{t,t',s} \leq 0;\ \mathrm{SF}_{l,s} + \underline{D}_{l,s}^{t'} + \underline{T}_{l,s,q}^{t} - \underline{T}_{l,s,q}^{t} \leq \mathrm{SF}_{l,q};\ \ 2\mathrm{SF}_{l,q} + \underline{D}_{l,q}^{t'} < \underline{D}_{l,q}^{t};\ \ \underline{D}_{l,q}^{j} \leq 0 \tag{34}$$

$$x_{t,t',s} \leq 0;\ \mathrm{SF}_{l,s} + \underline{D}_{l,s}^{t'} + \underline{T}_{l,s,q}^{t'} - \underline{T}_{l,s,q}^{t} > \mathrm{SF}_{l,q};\ \ 2\mathrm{SF}_{l,q} + \underline{D}_{l,q}^{t'} \geq \underline{D}_{l,q}^{t};\ \ 0 \leq \underline{D}_{l,q}^{t'} + \mathrm{SF}_{l,q} \tag{35}$$

$$x_{t,t',s} \leq 0;\ \mathrm{SF}_{l,s} + \underline{D}_{l,s}^{t'} + \underline{T}_{l,s,q}^{t'} - \underline{T}_{l,s,q}^{t} > \mathrm{SF}_{l,q};\ \ 2\mathrm{SF}_{l,q} + \underline{D}_{l,q}^{t'} < \underline{D}_{l,q}^{t};\ \ \underline{D}_{l,q}^{t} \leq \underline{D}_{l,q}^{t'} \tag{36}$$

As stated before, these rules are able to fix some binary variables. Thus, by the reduction of binary variables a solution in a shorter time can be obtained; however, there is no guarantee about its optimality and it can act as an upper-bound for the problem. In addition, if we consider the trains in a group, that they have same travel and dwell times, the phrase $\underline{T}_{l,i,s,s+1} - \underline{T}_{l,j,s,s+1}$ will relax and the mentioned conditions for each constraint will simplify. Furthermore, the dwell times for most of the stations in a network are equal and can be relaxed from some of the constraints. Only in the stations that there is a crossing point of two or more lines, dwell time is slightly different with other stations. In addition, the safety time, which is based on the geographical specifications of a station and the line, usually is equal for stations. Thus, the mentioned constraints are simplified for most trains and stations so that they are able to decrease the complexity of the problem.

## 4- Lower Bound

In this section, we present a Lagrangian Relaxation (LR) lower bound algorithm to estimate a powerful lower bound of the objective function. Lagrangian Relaxation algorithm is one of the most efficient algorithms, obtaining lower bound. In this algorithm complex constraints are relaxed from the set of constraints and with a penalty multiplier are added to the objective function. Selecting the relaxed constraints with non-zero integrality gap and updating the multiplier of relaxed constraint are of high importance [54]. Here, the overtaking constraints and capacity related constraints increase the complexity of the problem. The binary decision variables are the elements which increase the complexity of the constraints. Constraint (12) binds all the binary decision variables and it seems that it is the most difficult constraint in the current set. Thus, constraint (12), which integrates the events together, is the best candidate, satisfies the mentioned criteria, and should be selected so that the problem efficiently relaxes. Considering the selected constraint, the objective function of LR model is as follows.

$$\min z = \sum_{l=1}^{L}\sum_{t=1}^{T} c_{ltm} + \sum_{l \in L_E}\sum_{t \in T_l}\sum_{t' \in T_l; \, t \neq t'}\sum_{s \in S_l} u_{l,t,t',s}^{k}(x^{"}{}_{t,t',s} + x'{}_{t,t',s} + x_{t,t',s} + x_{t',t,s} - 1) \tag{37}$$

which, $u_{l,t,t',s}^{k}$ are the Lagrangian multipliers and the sub-gradient method updates them because of its efficacy [55, 56, 57]. Other constraints, except the constraint (12), are embedded in the LR model. The multipliers initially are interpreted as the price (marginal cost) of the relaxed constraint in a feasible solution. Because the constraint is in the equality form, the multipliers in a feasible solution always are zero. The multipliers are iteratively adjusted using the result of the model in a way that helps to improve the amount of the lower bound. Therefore, in each step amount of $x^{"}{}_{t,t',s} + x'{}_{t,t',s} + x_{t,t',s} + x_{t',t,s} - 1$ is calculated and called as $\gamma_{l,t,t',s}$, and it is used to update the multipliers. Equations (38), (39) and (40) demonstrate this process as sub-gradient method.

$$\gamma_{l,t,t',s} = x^{"}{}_{t,t',s} + x'{}_{t,t',s} + x_{t,t',s} + x_{t',t,s} - 1 \tag{38}$$

$$Step\_Size_k \leq \theta \frac{(UB - Objective_k)}{\sum_{l \in L_E}\sum_{t \in T_l}\sum_{t' \in T_l; \, t \neq t'}\sum_{s \in S_l}(\gamma^2{}_{l,t,t',s})} \tag{39}$$

$$u_{l,t,t',s}^{k+1} = u_{l,t,t',s}^{k} + \gamma_{l,t,t',s} \times Step\_Size_k \tag{40}$$

where $k$ is the iteration index used in the LR model, $UB$ is the objective function of a feasible solution, "$Objective_k$" is the amount of the objective function of the LR model in iteration $k$ and $\theta$ is a parameter by an initial value as 2. Equation (38) defines the amount that relaxed constraint is not satisfied. Equation (39) calculates the amount of the step-size parameter in iteration $k$ and equation (40) updates the amount of Lagrangian multipliers for the next iteration. Using these parameters, the LR iteratively updates the parameters to reach the optimal solution. The related procedure is shown in Figure 17.

[Figure 5 should be here]

At the first step, in order to define *UB*, which is used in equation (39), we find a feasible solution in which trains go through the stations with a lexicographical order and without any overtaking. The objective function of this solution that multiplies to 1.05 defines *UB*. Moreover, we get the objective function of the problem with relaxed binary variables. The objective function of this Relaxed Mixed Integer Problem (RMIP) will be used as a benchmark for calculating the improvement of the LR model. After the first iteration, the parameters and multipliers will update and the next iteration will start until one of these conditions meet:

- $\max\limits_{l, t, t', s} (u_{l,t,t',s}^{k+1} - u_{l,t,t',s}^{k}) < 0.005$
- Number of iterations exceeds 100.

Moreover, in each step if the objective function of LR model does not increase, the amount of $\theta$ will be halved, which helps it to increase the lower bound.

## 5- Experimental Results

In order to demonstrate the computational efficiency of the mathematical programming model and the proposed upper and lower bound rules, a series of numerical experiments are illustrated. The following case study is based on the real parameters of Tehran Metro line 5 between Tehran and Golshahr with 12 stations and 40 kilometers length, which its general view is magnified in Figure 18. The actual Tehran−Golshahr line is a double-tracked line, in which every 11-minutes trains run from Tehran to Golshahr and vice versa. There are some express trains that stop only in three stations and they have no dwell time in other stations. The travel time for the express trains is about 32 minutes and it is 52 minutes for other trains.

Each station has at least one platform in which loading and unloading of passengers occurs. Except for the first and last stations, the main line of each station has no platform for passengers to board or alight from trains. Trains can overtake each other in the sidings or in the stations.

There are different headways in the network in a day. The network headway in the peak-hours alters and the minimum amount of the departure time between two trains reaches 480 seconds. Regarding the complexity of timetabling in the peak-hours, we consider related parameters of the peak-hours to obtain a timetable.

The computational time of train timetabling algorithms can be affected by number of siding's and station's capacities, number of trains, number of stations and sidings and the variability of trains' parameters. In the following numerical experiments, we focus on the impact of the variability of station capacity, trains and number of trains, since these factors mainly influence the structure of a train timetable and the resulting number of possible solutions.

However, the real-world problem used in this study only offers two types of trains, and a fixed number of trains that creates a simple problem to solve. To allow a comprehensive and systematic assessment, we construct random instances to evaluate the performance of the proposed algorithm. In this way, we increased the capacity of some middle stations. The stations were chosen randomly and the number of facilities increased randomly up to three. As another factor for increasing the complexity of the problem, we increased the variability of the travel time for some trains in a way that for 50% of randomly selected trains we increased the travel time up to 100% for all stations.

In the following experiments, consciousness and authenticity of the model, the effectiveness of the heuristic upper bound rules and the quality of the proposed Lagrangian Relaxation method are investigated. In these investigations five problems with different number of trains have been solved with the proposed models. We have used GAMS workstation and CPLEX solver to solve the problems on a PC equipped with 3 GHz sixteen cores processor and 18 GB of RAM.

[Figure 6 should be here]

Figure 19 illustrates the resulting optimal timetable for nine trains for the case of Tehran- Golshahr line. As shown in the Atmosfer station, three trains 3, 4, and 5 simultaneously are in the station, in station Vardavard and Garmdareh train 3 overtakes trains 5 and 4, and the other general constraints of the problem are satisfied.

[Figure 7 should be here]

5-1  Performance of Upper-Bound Rules

In order to analyze the five introduced upper-bound rules, we solved each problem by GAMS using the CPLEX solver, as an efficient and exact solver, and the optimal solutions gathered as a benchmark. The optimal solutions are shown in Table 4. The quality of the rules is measured by percentage gap between the obtained upper-bound and the corresponding optimal value. In addition, improvement in CPU time is the second most important criterion in analyzing the obtained rules.

In order to analyze the effectiveness of upper-bound rules, the obtained constraints (constraints (23), (33), (34), (35) and (36)) are embedded in the problem. This problem is named as Upper-Bound Rule Problem (UBRP) in the results and hereafter. The objective function and CPU time of UBRP is also shown in Table 4.

[Table 1 should be here]

As shown in Table 4, in all discussed problems when five upper-bound rules are added to the set of constraints the objective function of the problem is equal to the optimal solution, which is obtained by the CPLEX solver. Thus, there is no optimality gap between solutions. This shows that no area of optimal solution space is removed by these rules. On the other hand, the results show that the CPU time is considerably reduced. This reduction in CPU time is the result of reduction in the number of branchings and also increasing the speed at which UBRP closes the opened nodes in the Branch and Bound tree. For more explanation on Branch and Bound see [52]. These are the benefits of the proposed rules which are described as follows.

I.  Reduction in the number of branchings on the nodes.

In order to clarify, consider a node that the CPLEX solver branches on to create the two new sub-problems. The branching is on the value of variable $x_{t,t',s}$ which its value is fixed in the UBRP. Therefore, in UBRP the branching on this variable does not occur. If the value of variable $x_{t,t',s}$ in CPLEX solver is integer, the node in the UBRP is also feasible and compared to the CPLEX solver, the tree size from that node halves. Also, the node is infeasible in UBRP if the value of $x_{t,t',s}$ is not integer and so the tree does not expand more. Thus, by fixing binary variables the number of branching decreases and as a result the CPU time decreases.

II.  Reduction of time required to fathom the opened nodes.

Fixing the binary variables decreases the required time to fathom the opened nodes. After a node opens, regarding the rules, the node is infeasible or some of their binary variables are known. Thus, in UBRP the nodes fathom in a shorter time compared to CPLEX. This process increases the speed of Branch and Bound algorithm and in less time more nodes can be investigated.

For more investigation, three diagrams of "number of opened nodes versus duality gap", "CPU time versus duality gap" and "CPU time versus number of opened nodes" are shown in Figure 20, Figure 21 and Figure 22 consecutively. The duality gap in these figures refers to the gap between the upper and lower bounds, which can be obtained at each step of Branch and Bound algorithm.

[Figure 8 should be here]

[Figure 9 should be here]

As shown in Figure 20, in all problems with optimality gap equal to zero, the number of opened nodes in UBRP is lower than the CPLEX solver (as benefit I). Also, as you can see in Figure 21 the duality in UBRP converges to zero quicker than in CPLEX. This is because in a shorter time more nodes opened and fathomed as benefit II which we have mentioned earlier. As mentioned in benefit II, fixing binary variables increases the speed at which the nodes fathom, i.e. a higher number of nodes fathom in an equal time. For more investigations we proposed a measurement criterion as ξ= *CPU time/ Opened nodes* that measures the speed at which the opened nodes close. This criterion measures the impact of the proposed rules on the node fathoming speed. The amount of ξ for each problem is shown in Table 5.

[Table 2 should be here]

As is shown in Table 5 the criterion ξ is reduced for all problems. In all problems in order to reach a given duality gap, UBRP opens fewer nodes in a shorter period of time (except the problem with 12 trains that is stopped in 400000 seconds and results of CPLEX is not known). The reason for this, is as described in benefit I and shows that the proposed rules reduce the number of branching.

In order to clarify the improvement, consider the problem with nine trains. In this problem in order to reach the duality gap about 0.5%, CPLEX opened 305000 nodes in 206 seconds, however the UBRP opened 160000 nodes (48% of opened nodes in CPLEX) in 107 seconds (is equal 51% of CPLEX CPU time). From another view point, in the problem with eleven trains UBRP after about 3775 seconds of solving process opened 3057644 nodes and the gap is about 0.99% versus 6008691 opened nodes (is equal 51% of UBRP) in 52923 seconds (is equal 7.1% ,i.e. 1 to 14 of CPLEX) and 0.98% gap in CPLEX.

On the other hand, in an equal period of time UBRP opens more nodes because of improvement in ξ. Figure 23 shows "CPU time versus number of opened nodes". As shown, in all problems except the problem with eight trains, in an equal period of time UBRP opens more nodes and in that period of time creates an smaller duality gap. E.g. in the problem with ten trains after about 1290 seconds of solving process UBRP opened 719791 nodes and the gap is about 0.74% versus 151211 opened nodes (is equal 21% of UBRP) and 1.13% gap in CPLEX (1.52 times greater than UBRP).

In sum, we can conclude that UBRP opens more nodes in an equal period of time while achieving a better duality gap and also required number of the opened nodes to get the optimal solution is lower than CPLEX.

[Figure 10 should be here]

Last but not least, the number of discrete variables in CPLEX and UBRP is shown Figure 23.

[Figure 11 should be here]

Figure 23 also shows the trend of discrete variables in UBRP and CPLEX. It is obvious that the two trends go ahead while the difference between the numbers of discrete variables increases. In other words, by increasing the number of trains the number of relaxed discrete variables increases and that engenders the reduction of CPU time. Regarding the figure, the improvements in CPU time can be explained more clearly.

### 5-2  Lower Bound Results

In order to investigate the results of Lagrangian Relaxation algorithm for a lower-bound, we used the same problem, which is described in last section. The measurement criteria are the duality gap between optimal value and the lower bound and the CPU time. In order to obtain *UB*, a feasible solution is obtained as the mentioned procedure. Our algorithm starts with result of the relaxed mixed integer programing (RMIP) model as a base for LR comparison and tries to improve it. Maximum number of the iterations in our investigation set to 25. In our investigations, there is only a little amount of improvement after iteration 15 and the algorithm tries to prove the optimality of the solution.

Results of LR model in Table 6 show that it has achieved the optimal value of problem and their value are the same as the values that have been reported by CPLEX in Table 4. The achievement stipulates the idea that the relaxed constraint has been correctly selected.

[Table 3 should be here]

To investigate further the quality of the algorithm, the trend of the lower-bound in iterations of the LR model for problems on 6 and 7 trains is shown in Figure 24.

[Figure 12 should be here]

Considering the fact that the selected relaxed constraint is in the form of equation and bind all the defined events together, it seems to be a very powerful constraint. As Figure 24 shows, the LR algorithm iteratively improves the lower bound until it gets the optimal solution. In addition, because the relaxed constraint is in an equality form there is no fluctuation in the sign of $\gamma_{l,t,t',s}$; which, results in a straightforward trend of improvement without fluctuation.

Nonetheless, the CPU time as another quality measurement criterion has not improved. In order to improve the CPU time, we defined another model using constraint (16), which is a special ordered set (SOS) of constraints and considers the platforms capacity. A same experiment analysis was performed on that model. The result showed a fluctuation without any improving trend. Also, as another model we relaxed constraints (2) and (3). The idea was from the problem of Seoul metropolitan Railway network that [37] relaxed the constraint that connects the stations with each other. The result of their analysis was hopeful and encouraging; nevertheless, the model that the connecting constraints - constraints (2) and (3)- relaxed does not result in an admirable solution in our case. The result again showed a fluctuation without any improved trend.

**Conclusion**

In this study a new MIP model for the railway scheduling problem is proposed that utilizes events of the railway networks. The general specifications and restrictions of railway networks are considered and also the station capacity constraints are embedded in the timetabling problem simultaneously. The model is NP-hard and is reduced to a flow-shop problem with block, machine assignment restriction and start time restrictions. Therefore, solving large scale problems in an acceptable amount of time is not possible.

In order to reduce CPU time some upper-bound rules based on analysis of model's parameters have been presented. These rules consider the relation among parameters of the model and with a grasp view try to eliminate some events that rationally are not wise to occur. The results has been added as five constraints to the master problem and experimental results testify that CPU time reduces up to 94% and the optimal solution also has been achieved. Furthermore, the analysis of the opened nodes shows that the rules can fathom more nodes in a shorter time, which describes its efficiency. These rules can be used on any other models by mapping the related variables.

By relaxing the constraint that ensures only one of the events occurs each time, a Lagrangian Relaxation algorithm is proposed. In order to update the Lagrangian multipliers, we used the step size method. Numerical analysis admitted accuracy and efficiency of the proposed algorithm, in a way that in all of the samples the optimal solution were achieve after about 20 iterations. However, the solving time of LR model is higher than the main problem, it creates some hopes to improve the processing time by focusing more on the proposed Lagrangian Relaxation algorithm.

Our on-going research has focused in heuristic methods in order to find upper bounds. With regard to the new definition of decision variables, further research can be done on implementation of a decomposition method as Dantzig-Wolfe algorithm, Benders decomposition or Branch and Price. Also, integrating of the proposed rule in this article by a pricing algorithm will be worthwhile or using them in Bisection algorithm.

**Acknowledgements**

The authors have benefited from Tehran Metro Agency's motivation and helpful comments in the model. This article has benefited from the advice and encouragement of Prof. N. Salmasi. Also, we wish to thank the authors of all the cited papers who motivated our work in this area.

# *List of Figures*

# *List of Tables*

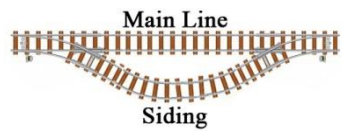

**Figure 13**

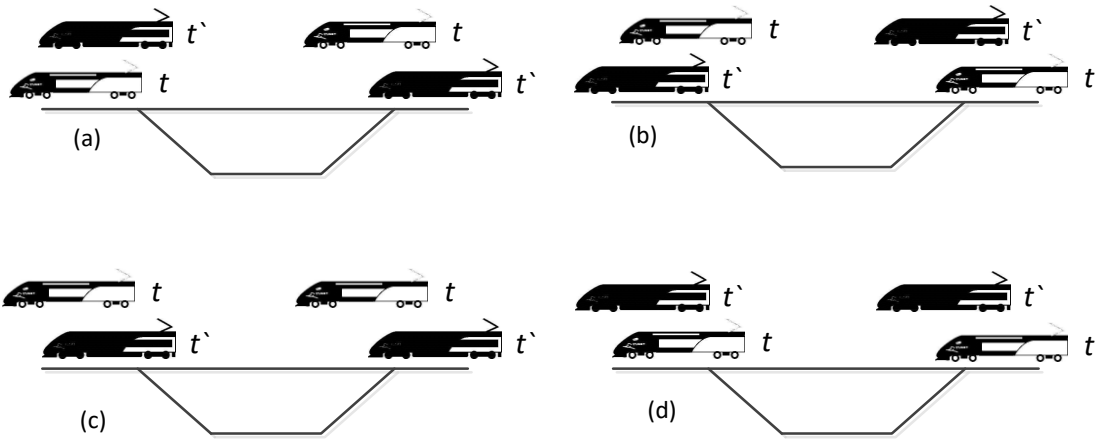

**Figure 14**

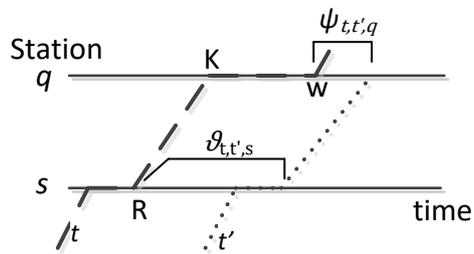

**Figure 15**

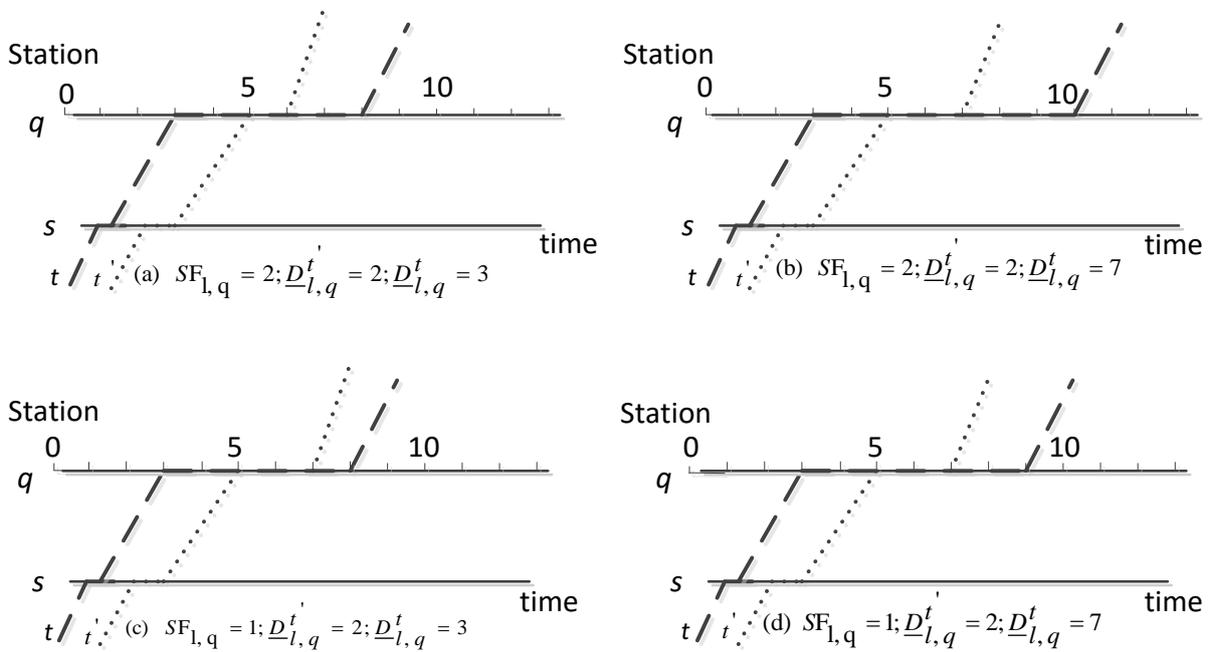

**Figure 16**

**Figure 17**

**Figure 18**

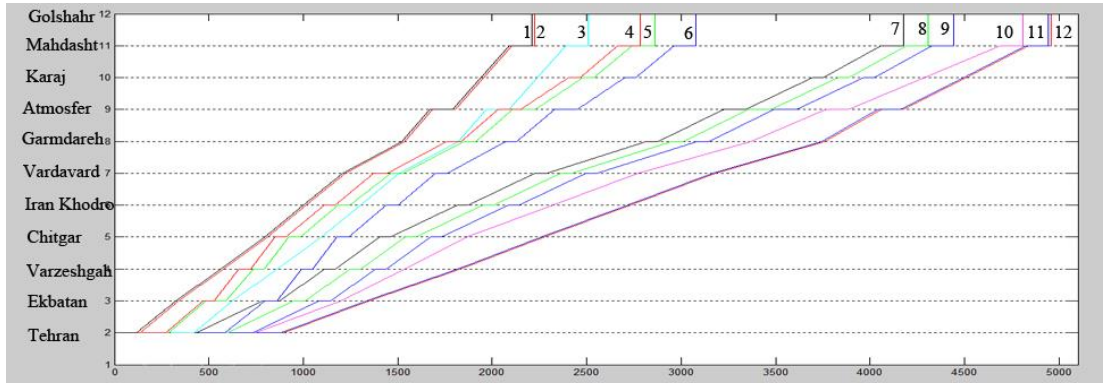

**Figure 19**

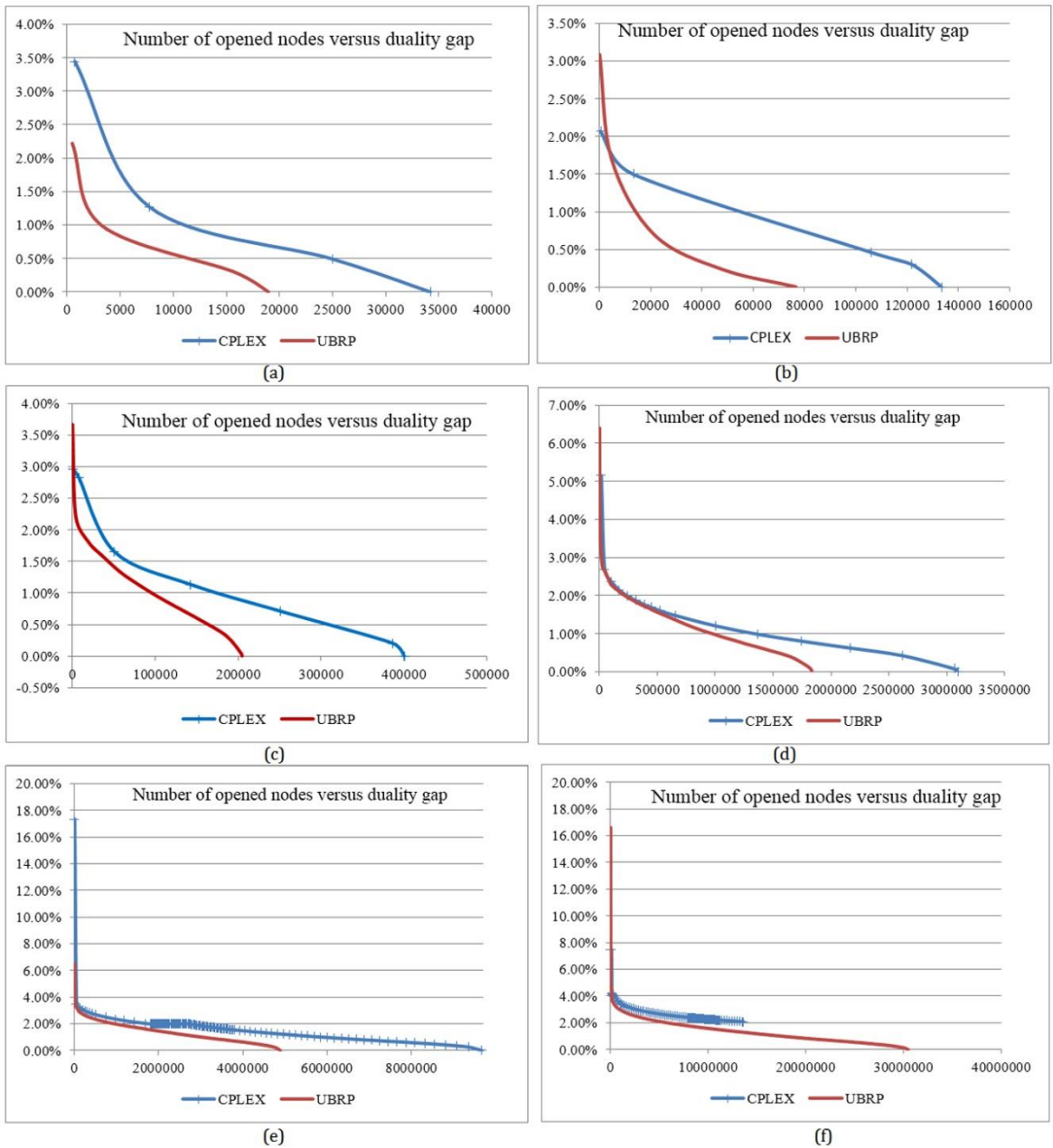

**Figure 20**

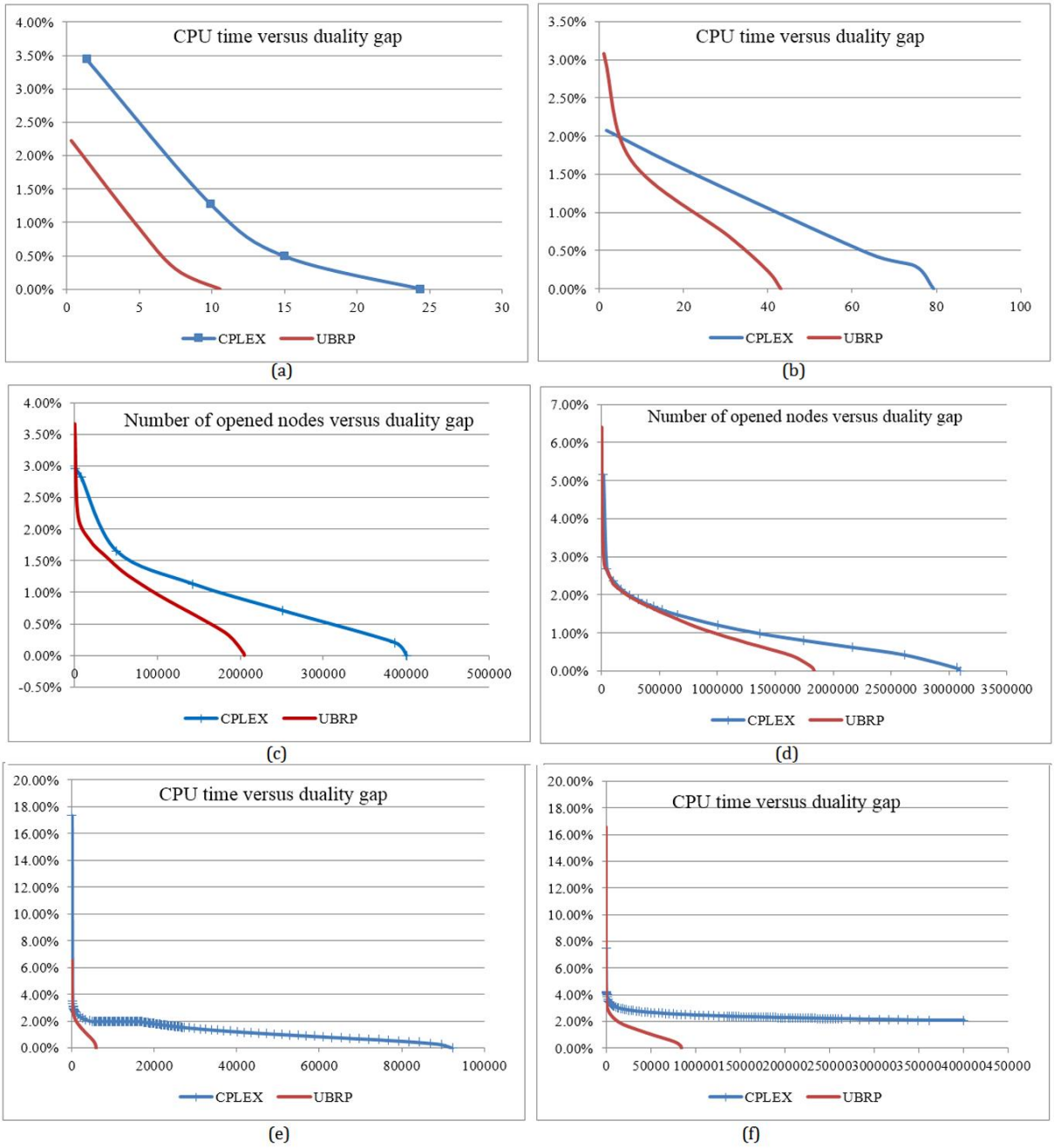

**Figure 21**

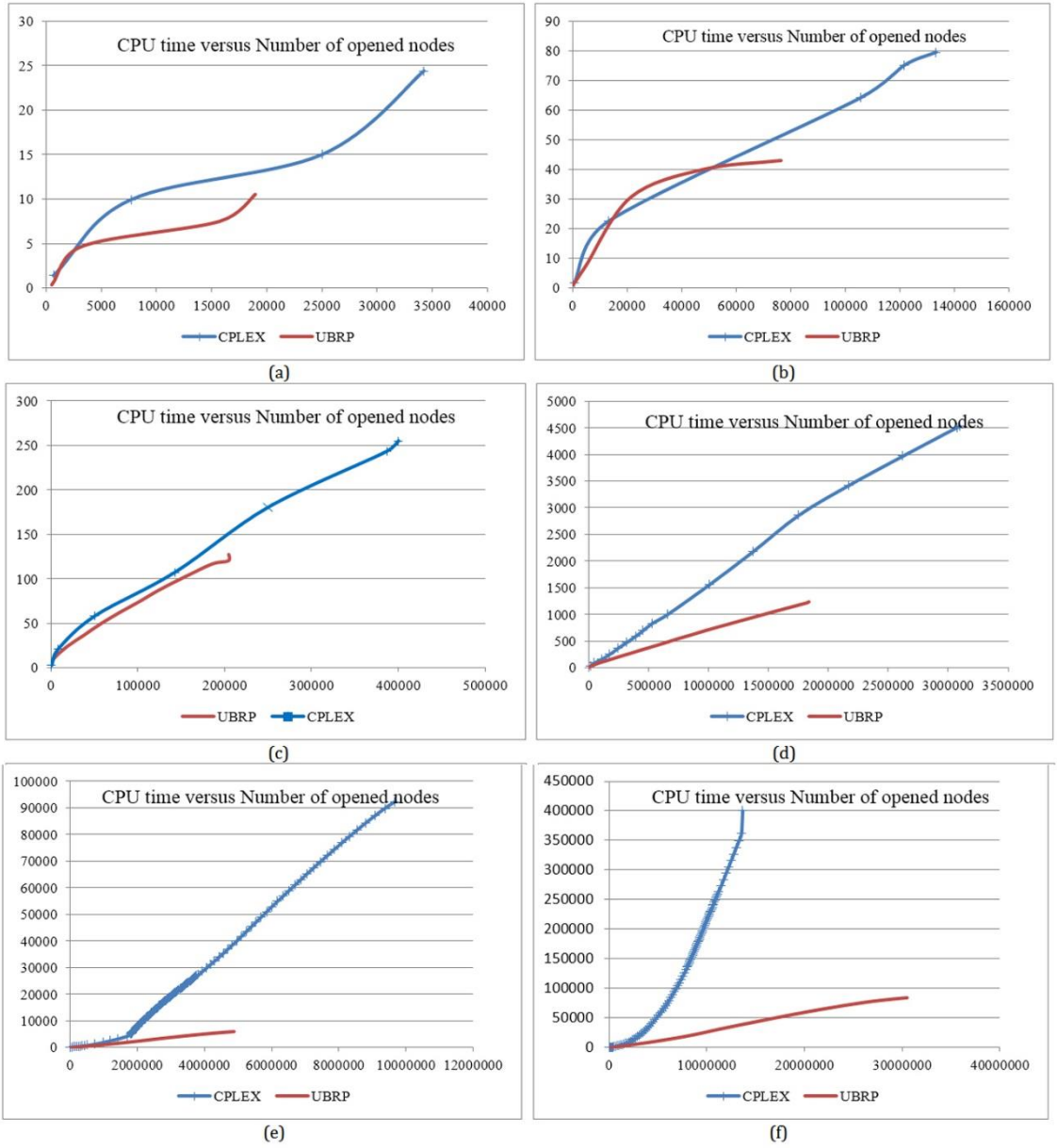

Figure 22

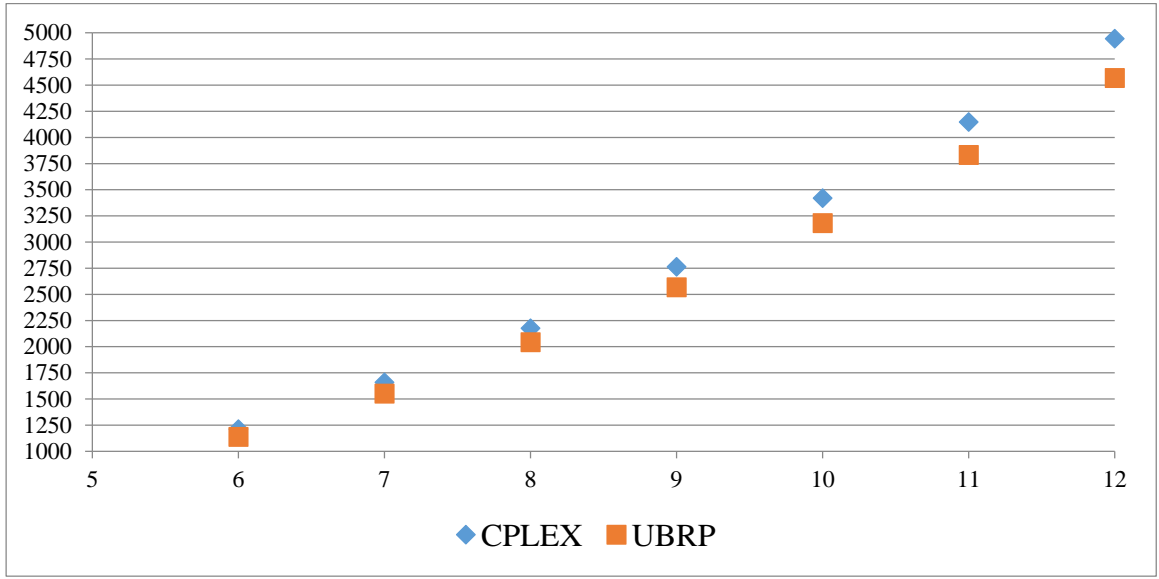

**Figure 23**

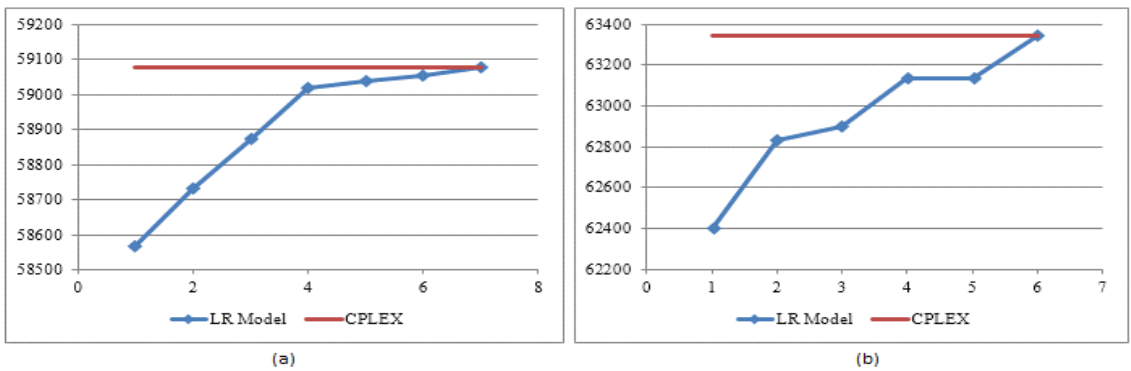

**Figure 24**

**Table 4**

| Number of trains | | 6 | 7 | 8 | 9 | 10 | 11 | 12 |
|---|---|---|---|---|---|---|---|---|
| Optimal solution | Objective | 63343 | 67683 | 70399 | 75133 | 78480 | 83410 | 86845 |
| (CPLEX) | CPU Time | 4.02 | 24.38 | 79.3 | 254.75 | 4526.3 | 92314 | 400008 |
| UBRP | Objective | 63343 | 67683 | 70399 | 75133 | 78480 | 83410 | 86845 |
| | CPU Time | 2.25 | 9.58 | 43.23 | 127.3 | 1234 | 5796.7 | 84276 |
| Optimality gap | | 0 | 0 | 0 | 0 | 0 | 0 | - |
| Improvement in CPU Time | | 44% | 61% | 45% | 50% | 73% | 94% | 79% |

**Table 5**

| Number of trains | | 6 | 7 | 8 | 9 | 10 | 11 | 12 |
|---|---|---|---|---|---|---|---|---|
| CPLEX | CPU Time (second) | 4.02 | 24.38 | 79.3 | 254.75 | 4526.4 | 92314 | 400084 |
| | Number of Opened Nodes | 13864 | 34226 | 133384 | 400422 | 3097165 | 9672652 | 13630809 |
| | CPU Time/ Opened Nodes ($\xi$) | 0.0019 | 0.0015 | 0.0016 | 0.0166 | 0.0481 | 0.0095 | 0.0294 |
| | Number of Discrete Variables | 1212 | 1659 | 2176 | 2763 | 3420 | 4147 | 4944 |
| UBRP | CPU Time (second) | 2.25 | 9.58 | 43.23 | 127.3 | 1234 | 5796.75 | 84276.3 |
| | Number of Opened Nodes | 6356 | 18956 | 76532 | 204396 | 1835564 | 4890725 | 30505371 |
| | CPU Time/ Opened Nodes ($\xi$) | 0.0004 | 0.0005 | 0.0006 | 0.0006 | 0.0007 | 0.0012 | 0.0028 |
| | Number of Discrete Variables | 1138 | 1552 | 2044 | 2568 | 3181 | 3833 | 4568 |
| Improvement percentage in $\xi$ | | 82% | 67% | 64% | 96% | 99% | 88% | 91% |

**Table 6**

| Number of trains | | 5 | 6 | 7 |
|---|---|---|---|---|
| LR Model | Best Bound | 59078 | 63343 | 67683 |
| | Number of Iteration | 6 | 7 | 22 |
| | CPU Time | >3.35 | >4.02 | >24.38 |

**Author's Biographical Sketch**

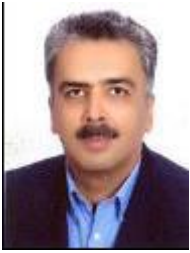

Dr. Kourosh Eshghi is a professor of Industrial Engineering Department at Sharif University of Technology, Tehran, Iran. He received his Ph.D. in the field of Operations Research from University of Toronto in 1997. He was the Head of IE Dept. for 7 years. His field of  interests include graph theory, integer programming and facility location. He is the author of 4 books and over 80 international journal papers.

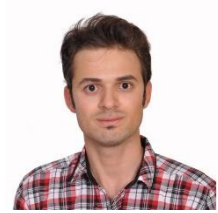

Afshin Oroojlooy jadid: is a PhD Student of Industrial Engineering in the department of Industrial Engineering at Lehigh University, Bethlehem, PA, USA. He received his MSc from the Sharif University of Technology in the Industrial Engineering and his BSc in industrial Engineering from Isfahan University of Technology. His main areas of research interests include Reinforcement Learning, Supply Chain Management, and Mixed Integer Programming.